\documentclass[12pt]{article}
\usepackage{amssymb}
\newcommand{\sect}[1]{\setcounter{equation}{0}\section{#1}}

\newcommand{\be}{\begin{equation}}
\newcommand{\ee}{\end{equation}}
\newcommand{\bea}{\begin{eqnarray}}
\newcommand{\eea}{\end{eqnarray}}
\newcommand{\beano}{\begin{eqnarray*}}
\newcommand{\eeano}{\end{eqnarray*}}
\newcommand{\nonu}{\nonumber \\}

\newcommand{\hs}[1]{\hspace{#1 mm}}

\newcommand{\vph}{\varphi}

\newcommand{\ca}{\mbox{$\cal{A}$}}

\newcommand{\cf}{\mbox{${\cal F}$}}
\newcommand{{\cg}}{\mbox{$\cal{G}$}}

\newcommand{\cR}{\mbox{$\cal{R}$}}

\newcommand{\cs}{\mbox{$\cal{S}$}}
\newcommand{\ct}{\mbox{$\cal{T}$}}
\newcommand{\cu}{\mbox{${\cal U}$}}
\newcommand{\cv}{\mbox{${\cal V}$}}

\newcommand{\prt}{\partial}

\newcommand{\wh}[1]{\widehat{#1}}
\newcommand{\wt}[1]{\widetilde{#1}}
\newcommand{\mb}[1]{\hs{4}\mbox{#1}\hs{4}}

\newcommand{\rmat}{{$R$-matrix}}
\newtheorem{theo}{Theorem}[section]
\newtheorem{prop}[theo]{Property}
\newtheorem{coro}[theo]{Corollary}
\newtheorem{defi}[theo]{Definition}

\newtheorem{lem}[theo]{Lemma}
\newcounter{rem}
\newtheorem{rmk}[rem]{Remark}
\newcommand{\prf}{\underline{Proof:}\ }
\newcommand{\finprf}{\null \hfill {\rule{5pt}{5pt}}\\[2.1ex]\indent}
\newcommand{\ie}{{\it i.e.}\ }
\newcommand{\CC}{\mbox{${\mathbb C}$}}

\newcommand{\ZZ}{\mbox{${\mathbb Z}$}}

\newcommand{\II}{\mbox{${\mathbb I}$}}

\begin{document}
\renewcommand{\thefootnote}{\fnsymbol{footnote}}
\newpage
\pagestyle{empty}
\setcounter{page}{0}

\newcommand{\LAP}{LAPTH}
\def\logo{{\bf {\huge LAPTH}}}

\centerline{\logo}

\vspace {.3cm}

\centerline{{\bf{\it\Large 
Laboratoire d'Annecy-le-Vieux de Physique Th\'eorique}}}

\centerline{\rule{12cm}{.42mm}}

\vspace{20mm}

\begin{center}

  {\LARGE  {\sffamily Vertex operators for quantum groups\\[.23ex]
  and application to integrable systems }}\\[1cm]

{\large E. Ragoucy\footnote{ragoucy@lapp.in2p3.fr}}\\[.21cm] 
 Laboratoire de Physique Th{\'e}orique \LAP\footnote{UMR 5108 
    du CNRS, associ{\'e}e {\`a} l'Universit{\'e} de Savoie.}\\[.242cm]
    LAPP, BP 110, F-74941  Annecy-le-Vieux Cedex, France. 
\end{center}
\vfill\vfill

\begin{abstract}
Starting with any $R$-matrix with spectral parameter, obeying the 
Yang-Baxter equation and a unitarity condition, we construct 
 the corresponding infinite dimensional quantum group $\cu_{R}$ in term  of a deformed 
oscillators algebra $\ca_{R}$. The 
realization we present is an infinite series, very similar to a
vertex operator. 

Then, considering the integrable hierarchy naturally associated to 
$\ca_{R}$, we show that $\cu_{R}$ provides its integrals of motion.
The construction can be applied to any infinite dimensional quantum 
group, e.g. Yangians or elliptic quantum groups.

Taking as an example the $R$-matrix of $Y(N)$, the Yangian based on $gl(N)$, 
we recover by this construction the nonlinear Schr\"odinger equation 
and its $Y(N)$ symmetry.
\end{abstract}

\vfill
\rightline{\tt mathQA/0108207}
\rightline{\LAP-859/01}
\rightline{July 01}

\newpage
\pagestyle{plain}
\setcounter{footnote}{0}

\markright{\today\dotfill DRAFT\dotfill Vertex operator\dotfill }

\sect{Introduction}
The aim of this paper is to present a general construction of 
infinite dimensional quantum 
groups as explicit integrals of motions of integrable systems. The 
construction relies only on the existence of an evaluated $R$-matrix (with 
spectral parameter) which obeys the unitarity condition. Thus, it can be 
applied to any infinite dimensional quantum group.

To the
$R$-matrix, one can associate a ZF algebra $\ca_{R}$ \cite{ZF}, 
which, in a Fock space representation provides the asymptotic states 
of the model. The quantum group is then constructed as an infinite 
series in the ZF generators, and shown to commute with the Hamiltonian 
of the hierarchy. Thus, it generates the integrals of motion of the 
hierarchy. Moreover, since there is a natural action of the quantum 
group on the $\ca_{R}$ generators, its action on asymptotic states of 
the system is easily deduced.

Taking as an example the $R$-matrix of $Y(N)$, the Yangian based on $gl(N)$, 
we recover by this construction the nonlinear Schr\"odinger equation 
and its $Y(N)$ symmetry \cite{MRSZ,MW}. 
It is thus very natural to believe that the 
other integrable systems known in the literature can be treated with 
the present approach.

\null

The paper is organized as follows. In the  section \ref{s:def}, we 
introduce the different definitions and properties we will need. {From} these notions, 
 we construct, in section \ref{s:prop}, a quantum group $\cu_{R}$
 {from} the deformed oscillator algebra $\ca_{R}$. We consider 
 in section \ref{s:hier} the hierarchy associated to $\ca_{R}$ and show 
that $\cu_{R}$ generates  integrals of motion.
Then, its Fock space 
 representation is studied in section \ref{s:fock}. Section \ref{s:ex} 
deals with three examples: the nonlinear Schr\"odinger equation with its 
Yangian symmetry (case of additive spectral parameter $R$-matrix), and 
$\cu_{q}(\wh{gl_{2}})$ and $\ca_{q,p}(gl_{2})$ (case of multiplicative $R$-matrix). Finally, 
we conclude in section \ref{s:concl}.

\sect{Definitions and first properties\label{s:def}}
\subsection{Z.F. algebras}
We start with an $R$-matrix satisfying the Yang-Baxter equation with 
spectral parameter:
\be
R_{12}(k_1,k_2)R_{13}(k_1,k_3)R_{23}(k_2,k_3)=
R_{23}(k_2,k_3)R_{13}(k_1,k_3)R_{12}(k_1,k_2)\label{YBE}
\ee
and the unitarity condition
\be
R_{12}(k_1,k_2)R_{21}(k_2,k_1)=\II\otimes \II \label{unitarity}
\ee
$R$ is an $N^{2}\times N^{2}$ matrix.
Here and below we will denote for briefness
\be
R_{12}\equiv R_{12}(k_1,k_2) \label{nota-R}
\ee
but let us stress that the $R$-matrix we consider are defined with spectral 
parameter. Note also that both the usual additive and multiplicative cases 
for the $R$-matrix, where $R(k_1,k_2)$ stands for $R(k_1-k_2)$ and 
$R(k_1/k_2)$ respectively, are included in our formalism.

\begin{defi}[ZF algebra $\ca_{R}$]\hfill\\
To each {\rmat} obeying (\ref{YBE}) and (\ref{unitarity}), one can 
associate a Zamolodchikov-Faddeev (ZF) algebra $\ca_{R}$ \cite{ZF}, with generators 
$a_{i}(k)$ and $a^\dag_{i}(k)$ ($i=1,\ldots,N$) and exchange relations:
\bea
a_{1} a_{2} &=& R_{21} a_{2} a_{1} \label{AN-1}\\
a_{1}^\dag a_{2}^\dag &=& a_{2}^\dag a_{1}^\dag R_{21} \label{AN-2}\\
a_{1} a_{2}^\dag &=& a_{2}^\dag R_{12} a_{1} +\delta_{12}
\label{AN-3}
\eea
We have used the notations
\beano
&&a_{1}=\sum_{i=1}^Na_{i}(k_{1})\, e_{i}\otimes\II
\mb{,} a_{2}=\sum_{i=1}^Na_{i}(k_{2})\, \II\otimes e_{i}\\
&&a^\dag_{1}=\sum_{i=1}^Na^\dag_{i}(k_{1})\, e^\dag_{i}\otimes\II
\mb{,} a^\dag_{2}=\sum_{i=1}^Na^\dag_{i}(k_{2})\, \II\otimes e^\dag_{i}\\
&& \delta_{12}=\delta(k_{1}-k_{2})\sum_{i=1}^N\, e_{i}\otimes e^\dag_{i} 
\mb{,} e^\dag_{i}=(0,\ldots,0,\stackrel{i}{1},0,\ldots,0)
\mb{,} e^\dag_{i}\cdot e_{j}=\delta_{ij} 
\eeano
where $\cdot$ stands for vectors scalar product.
\end{defi}
Let us remark that, in the same way the Yang-Baxter equation ensures 
the associativity of the product in $\ca_{R}$,
the unitarity condition can be interpreted as a 
consistency condition for the $\ca_{R}$ algebra. Indeed, starting with 
(\ref{AN-1}), exchanging the 
auxiliary spaces
$1\leftrightarrow 2$ and the spectral parameters 
$k_1\leftrightarrow k_2$, and multiplying by $(R_{12})^{-1}$ one gets
\be
a_{1} a_{2} = (R_{12})^{-1} a_{2} a_{1} 
\ee
Comparing this last relation with (\ref{AN-1}), we recover the unitarity 
condition.

Above and in the following, we will loosely write $a_{1}\in\ca_{R}$.
\begin{prop}[Adjoint anti-automorphism]\hfill\\
\label{adjoint}
\be
\mb{\it Let $\dag$ be the operation defined by}\left\{
\begin{array}{lcl}
        \ca_{R} & \rightarrow & \ca_{R}\\
a(k) & \mapsto & a^\dag(k)\\
a^\dag(k) &\mapsto & a(k)\\
R_{12}(k_{1},k_{2}) &\mapsto & R_{21}(k_{2},k_{1})
\end{array}\right.
\ee
and $(xy)^\dag=y^\dag x^\dag\
\forall\, x,y\in\ca_{R}$. Then $\dag$ is an automorphism of the $\ca_{R}$ 
algebra, 
and we can make the 
identifications $(a)^\dag\equiv a^\dag$ and $(a^\dag)^\dag\equiv a$.
\end{prop}
\prf
Direct calculation. For instance:
\beano
(a_{1}a_{2})^\dag &=& (a_{2})^\dag (a_{1})^\dag\\
&=& (a_{1})^\dag (a_{2})^\dag (R_{21})^\dag= (a_{1})^\dag 
(a_{2})^\dag R_{12}
\eeano
After the exchange $1\leftrightarrow2$, one recovers 
(\ref{AN-2}):
\be
(a_{1})^\dag (a_{2})^\dag=(a_{2})^\dag (a_{1})^\dag R_{21}
\ee
The other relations are obtained in the same way, once one remarks
$\big(\delta_{21}\big)^\dag=\delta_{12}$.
\finprf
\subsection{Vertex operators}
\begin{defi}[Vertex operators]\hfill\\
        \label{def.ov}
The vertex operators $T^{ij}(k)$, ($i,j=1\ldots,N$) 
associated to the algebra $\ca_{R}$ are defined by 
$T(k)\equiv T^{ij}(k)E_{ij}\in\ca_{R}\otimes \CC^{N^{2}}$ where
\be
T(k_{\infty}) = \II+\sum_{n=1}^\infty\, \frac{(-1)^n}{n!} 
a^\dag_{{n}\ldots {1}}\, T^{(n)}_{\infty1\ldots n}a_{1\ldots n}
\label{defT}
\ee
with 
\bea
 a^\dag_{{n}\ldots {1}} &=&\!  
a^\dag_{\alpha_{n}}(k_{n})\ldots a^\dag_{\alpha_{1}}(k_{1})\\
a_{1\ldots n} &=&\!  a_{\beta_{1}}(k_{1})\ldots a_{\beta_{n}}(k_{n}) \\
T^{(n)}_{\infty1\ldots n} &=&\!  T^{(n)}_{\infty,\alpha_1,\beta_{1},\ldots 
\alpha_n,\beta_{n}}(k_{\infty},k_1,\ldots,k_n)\in \left(\CC^{\otimes 
N^2}\right)^{\otimes (n+1)}\! (k_{\infty},k_1,\ldots,k_n)
\eea
In (\ref{defT}), there is an implicit summation on the indices
$\alpha_1,\beta_{1},\ldots,\alpha_n,\beta_{n}=1,\ldots,N$ and 
an integration over the 
spectral parameters $k_{1},\ldots, k_{n}$.

For convenience, 
$\infty$ labels the auxiliary space associated to $T(k_{\infty})$, 
and, as for the \rmat, we will note $T_{\infty}\equiv T_{\infty}(k_{\infty})$.
\end{defi}
Let us stress that, in the notation (\ref{defT}), the auxiliary 
spaces $1,\ldots,n$ are "internal" in the sense that the indices 
corresponding to these spaces are summed and define scalars, not 
matrices, in these spaces. It is only the indices 
corresponding to the "external" auxiliary space $\infty$ which refers 
to the matrix labeling for $T$.
For instance 
$a^\dag_{{1}}\, T^{(1)}_{\infty1}a_{1}$ stands for 
\beano
a^\dag_{{1}}\, T^{(1)}_{\infty1}a_{1} &=& 
\sum_{\alpha,\beta=1}^{N}\left(a^\dag_{{1}}\, 
T^{(1)}_{\infty1}a_{1}\right)_{\alpha,\beta}E_{\alpha,\beta}
\ =\ \sum_{\alpha,\beta=1}^{N}\left(\sum_{\gamma,\mu=1}^{N}
a^\dag_{\gamma}\, 
T^{(1)}_{\alpha,\beta;\gamma,\mu}a_{\mu}\right)E_{\alpha,\beta}
\eeano
so that we could have written $a^\dag_{{2}}\, T^{(1)}_{\infty 2}a_{2}$ 
as well: $1,\ldots,n$ are dummy 
space indices.

\null

\begin{rmk} \rm The series (\ref{defT}) is very similar to a normal ordered 
(in $a$ and $a^\dag$) exponential
\be
V(k_{\infty})\ =\ :\,\exp\left(-a^\dag M a\right)\, :
\ee
whence the denomination vertex operator used here to denote it.
\end{rmk}
\begin{prop}[$\cs_{n}$-covariance of the vertex operators]\hfill\\
\label{inv}
 The vertex operators coefficients $T^{(n)}_{\infty1\ldots n}$
 are covariant under the action of the permutation group $\cs_{n}$.
 
 More precisely, for $\sigma\in\cs_{n}$, one has:
 \be
 T^{(n)}_{\infty\sigma(1)\ldots \sigma(n)} = 
{\cR}^{1\ldots n}_{\sigma} T^{(n)}_{\infty 1\ldots n} ({\cR}^{1\ldots 
n}_{\sigma})^{-1}
\label{Sn-inv}
\ee
where ${\cR}^{1\ldots n}_{\sigma}$ is the product of $R$-matrices defined by 
$\displaystyle a_{\sigma(1)\ldots \sigma(n)} = 
{\cR}^{1\ldots n}_{\sigma} a_{1\ldots n}$.
\end{prop}
\prf
Starting {from} the term  
$X_{n}=a^\dag_{{n}\ldots {1}}\, T^{(n)}_{\infty1\ldots n}a_{1\ldots n}$
and relabeling the 
 auxiliary spaces $i\rightarrow \sigma(i)$ (and also the spectral 
 parameters), one gets
\be
X_{n}=a^\dag_{\sigma({n})\ldots \sigma({1})}
T^{(n)}_{\infty \sigma(1)\ldots \sigma(n)}a_{\sigma({1})\ldots \sigma({n})}
\ee
Then, {from} the exchange properties of the $a$'s and $a^\dag$'s
and the property \ref{adjoint}, one 
has:
\be
a^\dag_{\sigma({n})\ldots \sigma({1})}=a^\dag_{{n}\ldots {1}}
({\cR}^{1\ldots n}_{\sigma})^{-1}
\mb{and}
a_{\sigma({1})\ldots \sigma({n})}={\cR}^{1\ldots n}_{\sigma}a_{1\ldots n}
\ee
which leads to the formula (\ref{Sn-inv}).
\finprf
As an example, if $\sigma$ is just the transposition 
$i\leftrightarrow i+1$, one gets $\cR^{1\ldots n}_{\sigma}=R_{i,i+1}$ and the 
formula
\be
 T^{(n)}_{\infty1\ldots i-1,i+1,i,i+2\ldots n} = 
{R}_{i,i+1} T^{(n)}_{\infty 1\ldots n} {R}_{i+1,i}
\ee
\begin{prop}
The matrices $\cR^{1\ldots n}_{\sigma}$, $\sigma\in\cs_{n}$, defined by 
\be
a_{\sigma(1)\ldots \sigma(n)} = 
{\cR}^{1\ldots n}_{\sigma} a_{1\ldots n}
\ee
obey to 
\be
\cR^{\mu(1)\ldots \mu(n)}_{\sigma}\cR^{1\ldots n}_{\mu}=\cR^{1\ldots n}_{\sigma o \mu} 
\mb{so that}
(\cR^{1\ldots n}_{\sigma})^{-1}=\cR^{\sigma(1)\ldots \sigma(n)}_{\sigma^{-1}}
\label{Rsigmu}
\ee
{From} any matrix $M_{1\ldots n}\in (\CC^{N^{2}})^{\otimes n}$, one can 
construct a $\cs_{n}$-covariant one by
\be
\wt{M}_{1\ldots n}=\frac{1}{n!}\sum_{\sigma\in\cs_{n}}
   ({\cR}^{1\ldots n}_{\sigma})^{-1} M_{\sigma(1)\ldots \sigma(n)}{\cR}^{1\ldots n}_{\sigma}
   \label{symetrisation}
\ee
\end{prop}
\prf 
The first formula is proved by direct calculation:
\be
a_{\sigma \circ\mu(1)\ldots \sigma \circ\mu(n)} = 
{\cR}^{1\ldots n}_{\sigma \circ\mu} a_{1\ldots n}={\cR}^{\mu(1)\ldots \mu(n)}_{\sigma} 
a_{\mu(1)\ldots\mu(n)}={\cR}^{\mu(1)\ldots \mu(n)}_{\sigma}{\cR}^{1\ldots n}_{\mu}a_{1\ldots n}
\ee

Now, for the last formula, one has (for any $\mu\in\cs_{n}$):
\beano
\wt{M}_{\mu(1)\ldots\mu(n)}&=&\frac{1}{n!}\sum_{\sigma\in\cs_{n}}
   ({\cR}^{\mu(1)\ldots\mu(n)}_{\sigma})^{-1} M_{\sigma \circ\mu(1)\ldots \sigma \circ\mu(n)}
   {\cR}^{\mu(1)\ldots\mu(n)}_{\sigma}\\
   &=& \frac{1}{n!}\sum_{\sigma'\in\cs_{n}}
  ({\cR}^{\mu(1)\ldots\mu(n)}_{\sigma' \circ\mu^{-1}})^{-1} M_{\sigma'(1)\ldots \sigma'(n)}
   {\cR}^{\mu(1)\ldots\mu(n)}_{\sigma' \circ\mu^{-1}}
\eeano
where in the last expression, we have made the change of variable
$\sigma'=\sigma \circ\mu$. Now, using (\ref{Rsigmu}), one gets
${\cR}^{\mu(1)\ldots\mu(n)}_{\sigma' \circ\mu^{-1}}=
\cR^{1\ldots n}_{\sigma'}(\cR^{1\ldots n}_{\mu})^{-1}$, and
$({\cR}^{\mu(1)\ldots\mu(n)}_{\sigma' \circ\mu^{-1}})^{-1}=
\cR^{1\ldots n}_{\mu}(\cR^{1\ldots n}_{\sigma'})^{-1}$, so that  
$\wt{M}_{1\ldots n}$ is $\cs_{n}$-covariant.
\finprf

\begin{rmk}\rm 
    Strictly speaking, one can start with vertex operators which does 
    not obey the $\cs_{n}$-covariance (\ref{Sn-inv}), but the relevant 
    part in the vertex operator will be the covariant one, as given 
    by (\ref{symetrisation}). 
\end{rmk}
\subsection{Well-bred operators}
\begin{defi}[well-bred operators]\hfill\\
        An operator $L$ is said well-bred\footnote{We call these operators 
        "well-bred" because they act nicely (on $a$ and $a^\dag$).}
         (on $\ca_{R}$) when it
        acts on $a$ and $a^\dag$ as
\be
L_{1}a_{2}=R_{21}a_{2}L_{1} \mb{and} 
L_{1}a^\dag_{2}=a^\dag_{2}R_{12}L_{1}\label{eqLa}
\ee
\end{defi}
We give few properties of well-bred operators that will be useful in 
the following.
\begin{lem}\label{Ldag}
Let $L$ be a well-bred operator, then 
$L^\dag(k)L(k)$ is central in $\ca_{R}$.
\end{lem}
\prf
One applies the $\dag$ automorphism to the relations (\ref{eqLa}). We 
get:
\be
a^\dag_{2}L^\dag_{1}=L^\dag_{1}a^\dag_{2}R_{12} \mb{and} 
a_{2}L^\dag_{1}=L^\dag_{1}R_{21}a^\dag_{2}
\ee
Then a direct calculation shows that $L^\dag(k)L(k)$ commutes with 
$a$ and $a^\dag$. For instance
\[
L_{1}^\dag L_{1}a_{2}^\dag =L_{1}^\dag a_{2}^\dag R_{12}L_{1}
= a_{2}^\dag L_{1}^\dag L_{1}
\]
\finprf
\begin{lem}\label{lem-c}
        Let $L$ be a well-bred operator of $\ca_{R}$.
Then $c_{12}=L_{1}^{-1}L_{2}^{-1}R_{12}L_{1}L_{2}$ is central in 
$\ca_{R}$. It satisfies $c_{12}^{-1}=c_{21}$.
\end{lem}
\prf
Starting with (\ref{eqLa}), one gets
\be
L_{1}L_{2}a_{3}=R_{32}R_{31}a_{3}L_{1}L_{2}
\ee
which can be rewritten (after exchange $1\leftrightarrow2$) as
\be
R_{31}R_{32}a_{3}=L_{2}L_{1}a_{3}L_{1}^{-1}L_{2}^{-1}
\ee
Then
\bea
R_{12}L_{1}L_{2}a_{3} &=& R_{12}R_{32}R_{31}a_{3}L_{1}L_{2}\ =\ 
R_{31}R_{32}R_{12}a_{3}L_{1}L_{2}\ =\ 
R_{31}R_{32}a_{3}R_{12}L_{1}L_{2}\nonu
&=& L_{2}L_{1}a_{3}L_{1}^{-1}L_{2}^{-1}R_{12}L_{1}L_{2}
\eea
So that, multiplying by $L_{1}^{-1}L_{2}^{-1}$, we obtain
\be
c_{12}\, a_{3}=a_{3}\, c_{12}
\ee
Performing a similar calculation with $a_{3}^\dag$, we get 
$c_{12}\,a^\dag_{3}=a^\dag_{3}\,c_{12}$.

The last equation is a direct consequence of the unitarity condition.
\finprf

\sect{Construction of well-bred vertex operators\label{s:prop}}

We first give a characterization of well-bred vertex operators:
\begin{lem}\label{lem.eq}
The vertex operators $T$ is well-bred
if and only if $T^{(n)}_{\infty1\ldots n}$ obeys
\bea
&&\displaystyle{T^{(1)}_{\infty0}=\II-R_{\infty0}}\mb{and for } n\geq1:
\label{eqT}\\
&&\displaystyle{(n+1)\left\{T^{(n)}_{\infty1\ldots n}-
\left(\cR_{0,n}\right)^{-1}
R_{\infty0}T^{(n)}_{\infty1\ldots n}\cR_{0,n}\right\}\ =\ 
\sum_{i=1}^{n+1}
\left(\cR_{0,i-1}\right)^{-1}T^{(n+1)}_{\infty1\ldots n\vert i}\,\cR_{0,i-1}} 
\nonumber
\eea
where we have introduced
\[
\cR_{0,n}=\prod_{a=1}^{\longleftarrow\atop n}\,
R_{0a}\mb{;} T^{(n+1)}_{\infty1\ldots n\vert i}=
T^{(n+1)}_{\infty1\ldots i-1,0,i\dots n}\ (i\leq n)\mb{and}
T^{(n+1)}_{\infty1\ldots n\vert n+1}=
T^{(n+1)}_{\infty1\dots n0}
\]
\end{lem}
\prf
We prove the property by a direct calculation. We note 
$\wh{T}_{\infty}=T_{\infty}-\II$:
\beano
a_{0}\wh{T}_{\infty} &=& \sum_{n=1}^\infty\frac{(-1)^n}{n!}\, 
(a^\dag_{n}R_{0n}a_{0}+\delta_{0n})a^\dag_{n-1\ldots1}T^{(n)}_{\infty1\ldots n}
a_{1\ldots n}\\
&=& \sum_{n=1}^\infty\frac{(-1)^n}{n!}\Big\{\, a^\dag_{{n}\ldots {1}}R_{0n}\cdots 
R_{01}a_{0}\, 
T^{(n)}_{\infty1\ldots n}a_{1\ldots n}\Big.\\
&&\hs{7} +\Big.
\sum_{i=1}^n a^\dag_{{n}\ldots {i+1}}a^\dag_{i-1\ldots 1}R_{0n}\cdots 
R_{0i+1}\delta_{0i}T^{(n)}_{\infty1\ldots n}a_{1\ldots n}\Big\}
\eeano
Using 
\be
\delta_{0i}T^{(n)}_{\infty1\ldots n}a_{i}=T^{(n)}_{\infty1\ldots 
i-1,0,i+1\ldots n}a_{0}
\ee
and after a relabeling $j\rightarrow j-1$ for $j\geq i+1$, one gets
\beano
a_{0}\wh{T}_{\infty} &=& \sum_{n=1}^\infty\frac{(-1)^n}{n!}\, 
a^\dag_{{n}\ldots {1}}R_{0n}\cdots 
R_{01}a_{0}\, T^{(n)}_{\infty1\ldots n}a_{1\ldots n}- 
T^{(1)}_{\infty0}a_{0}\\
&&+\sum_{n=2}^\infty\frac{(-1)^n}{n!}\sum_{i=1}^{n} 
a^\dag_{{n-1}\ldots 1}R_{0n-1}\cdots 
R_{0i}T^{(n)}_{\infty1\ldots0i\ldots n-1}a_{1\ldots i-1}a_{0}a_{i\ldots n-1}
\eeano
with, as a notation: 
\[
 i=n:\ R_{0n-1}\cdots 
R_{0i}\equiv1, \ T^{(n)}_{\infty1\ldots0i\ldots n-1}\equiv 
T^{(n)}_{\infty1\ldots n-1,0}\mb{and}a_{1\ldots i-1}a_{0}a_{i\ldots 
n-1}\equiv a_{1\ldots n-1}a_{0}
\]
Rewriting
\be
a_{1\ldots i-1}a_{0}a_{i\ldots n-1}=
R_{0i-1}\cdots R_{01}a_{0}a_{1\ldots n-1}
\ee
and relabeling $n\rightarrow n-1$ in the second summation we are led 
to
\bea
a_{0}\wh{T}_{\infty} &=& -T^{(1)}_{\infty0}a_{0}+
\sum_{n=1}^\infty\frac{(-1)^n}{n!}\Big\{ a^\dag_{{n}\ldots {1}}R_{0n}\cdots 
R_{01}a_{0}\, T^{(n)}_{\infty1\ldots n}a_{1\ldots n}+ 
\Big.\\
&&\hs{23}-\frac{1}{n+1}\Big.\sum_{i=1}^{n+1} a^\dag_{{n}\ldots 1}R_{0n}\cdots 
R_{0i}T^{(n+1)}_{\infty1\ldots n\vert i}R_{0i-1}\cdots R_{01}
a_{0}a_{1\ldots n}
\Big\}
\nonumber
\eea
that is 
\bea
a_{0}\wh{T}_{\infty} &=& -T^{(1)}_{\infty0}a_{0}+
\sum_{n=1}^\infty\frac{(-1)^n}{n!}a^\dag_{{n}\ldots {1}}\Big\{ \cR_{0n} 
T^{(n)}_{\infty1\ldots n}+ 
\Big.\label{lhsaT}\\
&&\hs{23}-\frac{1}{n+1}\Big.\sum_{i=1}^{n+1}\cR_{0n}\cR_{0,i-1}^{-1}
T^{(n+1)}_{\infty1\ldots n\vert i}\cR_{0,i-1}
\Big\}a_{0}a_{1\ldots n}
\nonumber
\eea

On the other hand, one computes
\be
R_{\infty0}\wh{T}_{\infty}a_{0} = \sum_{n=1}^\infty\frac{(-1)^n}{n!} 
a^\dag_{{n}\ldots {1}}\, R_{\infty0}T^{(n)}_{\infty1\ldots n}\cR_{0n}
a_{0}a_{1\ldots n}\label{rhsTa}
\ee 
Finally, equaling (\ref{lhsaT}) and (\ref{rhsTa}), we get
the equations (\ref{eqT}), after left-multiplication by $\cR_{0n}^{-1}$.

A similar calculation on 
$T_{\infty}a^\dag_{0}=a^\dag_{0}R_{\infty0}T_{\infty}$
leads to the same equation.
\finprf

\begin{rmk}\rm  If one defines $\cR_{00}=\II$ (and $T^{(0)}_{\infty}=\II$ as 
given by (\ref{defT})), the equation $T^{(1)}_{\infty0}=\II-R_{\infty0}$ 
just corresponds to $n=0$ in (\ref{eqT}).
\end{rmk}

\begin{prop}[Central generators of $\ca_{R}$]\hfill\\
\label{centre}
The only central generators of $\ca_{R}$ are constants.
\end{prop}
\prf
Let $c$ be a central generator of $\ca_{R}$. Since it commutes with 
$a$ and $a^\dag$, it also commutes with the number operator
$H_{0}=\int dk\, a^\dag(k)a(k)$ (see section \ref{s:hier}). It is thus of the 
form
\be
c=c^{(0)}+\sum_{n=1}^\infty \frac{(-1)^n}{n!}\, a^\dag_{{n}\ldots {1}}\, 
c^{(n)}_{1\ldots n}a_{1\ldots n}
\ee
Demanding $c\,a_{0}=a_{0}c$ leads to equations on the elements
$c^{(n)}_{1\ldots n}$. These equations are computed in the same way 
one computes the equations for the $T^{(n)}$'s. Indeed, one deduces 
the equations on $c^{(n)}$ by formally  replacing $R_{0\infty}$ by 
$\II$ in the equations (\ref{eqT}). We get the relations
\be\begin{array}{l}
c^{(1)}_{1}=0\\
\displaystyle
(n+1)\left\{c^{(n)}_{1\ldots n}-\left(\cR_{0,n}\right)^{-1}
c^{(n)}_{1\ldots n}\cR_{0,n}\right\}\ =\ \sum_{i=1}^{n+1}
\cR_{0,i-1}^{-1}
c^{(n+1)}_{1\ldots n\vert i}\,\cR_{0,i-1}\mb{for} n\geq1
\end{array}
\label{eq.centre}
\ee
We prove by induction that $c^{(n)}=0$. The case $n=1$ is a direct 
consequence of the equations. Let us suppose that $c^{(p)}=0$ for 
$p\leq n$. Writing the equation (\ref{eq.centre}) at level $n$, and 
using the induction, we have
\be
\sum_{i=1}^{n+1} \cR_{0,i-1}^{-1}
c^{(n+1)}_{1\ldots n\vert i}\,\cR_{0,i-1}=0
\ee
Using the invariance property \ref{inv}, we can rewrite each term of the 
sum has
\be
c^{(n+1)}_{1\ldots n\vert i}=c^{(n+1)}_{1\ldots i-1,0,i,\ldots n}
=\cR_{0,i-1}c^{(n+1)}_{01\ldots n}\cR_{0,i-1}^{-1}
\ee
Thus, the equation is equivalent to $(n+1)\, c^{(n+1)}_{01\ldots n}=0$ 
and the induction is proven.
\finprf
\begin{theo}\label{theoTa}
        The vertex operators $T$ is well-bred
if and only if $T^{(n)}_{\infty1\ldots n}$ is
        defined by the 
        following inductive expressions:
\bea
T^{(1)}_{\infty0} &=& \II-R_{\infty0}\\
T^{(n+1)}_{\infty01\ldots n} &=& \frac{1}{n+1}
\sum_{i=0}^{n} (\cR^{01\ldots n}_{p_{i}})^{-1}\,
T^{(n)}_{\infty2\ldots i,0,i+1,..,n}\,\cR^{01\ldots n}_{p_{i}}\nonu
&&-\frac{1}{(n+1)!}
\sum_{\sigma\in\cs_{n+1}}
(\cR^{01\ldots 
n}_{p_{n}\circ\sigma})^{-1}\,R_{\infty\sigma(0)}T^{(n)}_{\infty\sigma(1)\ldots 
\sigma(n)} \,
\cR^{01\ldots n}_{p_{n}\circ\sigma}
\label{defTn}
\eea
where $T^{(n)}_{\infty2\ldots i,0,i+1,..,n}$ for $i=0$ stands for 
$T^{(n)}_{\infty1\ldots ,n}$.
 $p_{j}\in\cs_{n+1}$ is defined by 
\bea
p_{j}&:&
(0,1,\ldots,j-1,j,j+1,\ldots,n)\to(1,2,\ldots,j,0,j+1,\ldots,n),\ 
1\leq j\leq n \nonu
p_{0}&=&id\label{defpi}
\eea
We remind that $R_{ij}$ stands for $R_{ij}(k_{i},k_{j})$.
\end{theo}
\prf
We start with the lemma \ref{lem.eq} and show that $T$ obeys the above 
inductive expressions.
Remark that from the definition of $\cR_{0,i}$, one 
has 
\be
\cR_{0,i}a_{01\ldots n}=a_{1,2,\ldots,i,0,i+1\ldots n}\Rightarrow 
\cR_{0,i}=\cR^{01\ldots n}_{p_{i}}
\ee
where $p_{i}$ is defined by (\ref{defpi}).
We start from the equation (\ref{eqT}) and work with $\cs$-covariant 
matrices. Then, the right-hand side is equal to 
$(n+1)T^{(n+1)}_{\infty01\ldots n}$, while the left-hand side reads:
\be
\frac{1}{n!}\sum_{\sigma\in\cs_{n+1}}(\cR^{01\ldots n}_{\sigma})^{-1}
\left\{T^{(n)}_{\infty\sigma(1)\ldots\sigma(n)}-
(\cR^{\sigma(0)\ldots\sigma(n)}_{p_{n}})^{-1}R_{\infty\sigma(0)}
T^{(n)}_{\infty\sigma(1)\ldots\sigma(n)}\cR^{\sigma(0)\ldots\sigma(n)}_{p_{n}}\right\}
\cR^{01\ldots n}_{\sigma}
\ee

Now, we decompose $\cs_{n+1}$ with respect to $\cs_{n}$: any
$\sigma\in\cs_{n+1}$ is  of the form (for some $0\leq i\leq n$)
$\mu \circ p_{i}$ with\footnote{Strictly speaking, $\mu$ is still in 
$\cs_{n+1}$, but it obeys $\mu(0)=0$ so that its restriction to 
$[1,n]$ define an element of $\cs_{n}$.} $\mu\in\cs_{n}$ and $p_{i}$ 
defined in (\ref{defpi}).
Using the covariance of $T^{(n)}$, one gets for the first part of 
the right hand side:
\beano
rhs_1 &:=& \frac{1}{n!}\sum_{\sigma\in\cs_{n+1}}(\cR^{01\ldots n}_{\sigma})^{-1}
T^{(n)}_{\infty\sigma(1)\ldots\sigma(n)}\cR^{01\ldots n}_{\sigma}\\
&=&\frac{1}{n!}\sum_{i=0}^n
\sum_{\mu\in\cs_{n}}(\cR^{01\ldots n}_{p_{i}})^{-1}
(\cR^{p_{i}(0)p_{i}(1)\ldots p_{i}(n)}_{\mu})^{-1}
T^{(n)}_{\infty\mu(p_{i}(1))\ldots\mu(p_{i}(n))}
\cR^{p_{i}(0)p_{i}(1)\ldots p_{i}(n)}_{\mu}
\cR^{01\ldots n}_{p_{i}} \\
&=&\frac{1}{n!}\sum_{i=0}^n(\cR^{01\ldots n}_{p_{i}})^{-1}
\Big(\sum_{\mu\in\cs_{n}}
\ct^{(\mu)}_{\infty p_{i}(0)p_{i}(1)\ldots p_{i}(n)}\Big)
\cR^{01\ldots n}_{p_{i}} 
\eeano
with $\ct^{(\mu)}_{\infty 01\ldots n}=(\cR^{01\ldots n}_{\mu})^{-1}
T^{(n)}_{\infty\mu(1)\ldots\mu(n)}
\cR^{01\ldots n}_{\mu}$. Since $\mu(0)=0$, one has 
$\cR^{01\ldots n}_{\mu}=\cR^{1\ldots n}_{\mu}$.Then, using the 
$\cs_{n}$-covariance of $T^{(n)}$, one gets 
$\ct^{(\mu)}_{\infty 01\ldots n}=T^{(n)}_{\infty1\ldots n}$, 
$\forall\mu$, so that
\be
rhs_1=\sum_{i=0}^n(\cR^{01\ldots n}_{p_{i}})^{-1}
T^{(n)}_{\infty 2\ldots i,0,i+1\ldots n}
\cR^{01\ldots n}_{p_{i}} 
\ee
Finally, to get (\ref{defTn}), one remarks in the second sum of the r.h.s. 
that $\cR^{\sigma(0)\ldots\sigma(n)}_{p_{n}}
\cR^{01\ldots n}_{\sigma}=\cR^{01\ldots n}_{p_{n}\circ\sigma}$, due to 
(\ref{Rsigmu}).

The same calculation (done in reverse direction) also shows  that the 
inductive expressions obey the lemma \ref{lem.eq}.
\finprf
\begin{rmk} \rm Note that the inductive expression proves the unicity of the 
solution.
\end{rmk}
\begin{rmk} \rm The first terms in the series (\ref{defTn}) are
\beano
T^{(1)}_{\infty1} &=& \II-R_{\infty1} \\
T^{(2)}_{\infty12} &=& \II-R_{\infty2}+R_{\infty2}R_{\infty1}
-R_{21}R_{\infty1}R_{12}
\eeano
\end{rmk}
\begin{coro}\label{coroformT}
$\forall n\geq 0$, $T^{(n)}_{\infty 1\ldots n}$ is a non-vanishing 
polynomial of $R$-matrices. It has the form:
\be
T^{(n)}_{\infty 1\ldots n}=\II+\sum_{i=1}^n S^{(i)}_{\infty 1\ldots n}
\mb{with} 
S^{(i)}_{\infty 1\ldots n}=\sum_{\mu\in\cs_{n}} m_{\mu}
M_{\mu} \,
R_{\infty\mu(1)}\cdots R_{\infty\mu(i)}\,M_{\mu}^{-1}
\label{formdeT}
\ee
where $M_{\mu}$ are products of matrices $R_{ab}$ with $1\leq a,b\leq n$
and $m_{\mu}\in\ZZ$.
\end{coro}
\prf
We prove the corollary by induction. The explicit expressions given 
above prove that it is true for $n=0,1,2$. Now, suppose (\ref{formdeT}) is true up 
to $n$. Then, the equation (\ref{defTn}) shows that it is also true for 
$n+1$. Indeed, the two sums in (\ref{defTn}) have conjugation by 
$R$-matrices of type $M_{\mu}$. Moreover, only the first sum 
contributes to $\II$, and effectively leads to a coefficient 1, while 
the second sum increase the number of $R_{\infty a}$ 
($a=0,1,\ldots,n$) matrices by 1.
\finprf
\begin{rmk}\rm  The above formula shows that $T^{(n)}$ is invertible (as 
a series) for all $n$.
\end{rmk}

Using the theorem \ref{theoTa}, one can show
\begin{prop}\label{RTT}
        The well-bred vertex operators $T$ of theorem \ref{theoTa} obey FRT 
        relations:
        \be
        R_{12}T_{1}T_{2}=T_{2}T_{1}R_{12}, \ie 
        R_{12}(k_{1},k_{2})T_{1}(k_{1})T_{2}(k_{2})=
        T_{2}(k_{2})T_{1}(k_{1})R_{12}(k_{1},k_{2})
        \ee
        In other words, they generate an infinite dimensional quantum group with evaluated 
        $R$-matrix $R_{12}$. 
        In the following, we will denote this quantum 
        group $\cu_{R}$.
\end{prop}
\prf
We use the lemma \ref{lem-c} for $T$: 
$c_{12}=T_{1}^{-1}T_{2}^{-1}R_{12}T_{1}T_{2}$ is central in 
$\ca_{R}$ and such that 
\be
R_{12}T_{1}T_{2}=T_{2}T_{1}c_{12} \label{TTc}
\ee
$c_{12}$ being central, and due to the property \ref{centre},
it is a constant matrix $M_{12}$.
To 
identify the exact expression of $M_{12}$, we use the result of 
theorem \ref{theoTa}. Looking at (\ref{TTc}) as a 
series in the number of say $a$ operators and projecting on
 number 0, we get
$c_{12}=M_{12}=R_{12}$.
\finprf
\begin{rmk} \rm  Looking at the term linear in $a$, one gets 
\be
R_{12}(T^{(1)}_{13}+T^{(1)}_{23}-T^{(1)}_{13}\cdot T^{(1)}_{23}) =
(T^{(1)}_{13}+T^{(1)}_{23}-T^{(1)}_{13}\cdot T^{(1)}_{13}) c_{12}
\ee
Plugging into this equation the expressions of $T^{(1)}$ and 
$c_{12}$, one recovers the Yang-Baxter equation, which is indeed 
satisfied.
\end{rmk}

\begin{prop}\label{Tinv}
Let $T$ be the well-bred vertex operator of theorem \ref{theoTa}. Then, one 
has
\be
T^\dag(k)=T(k)^{-1}
\ee
\end{prop}
\prf
{From} the lemma \ref{Ldag}, one knows that $T^\dag(k)\,T(k)$ is 
central. This implies (using property \ref{centre})
that, $T^\dag(k) T(k)$ is a constant $N\times N$ matrix $M$. 
Looking at the term without $a$, one concludes 
that $M=\II_{N}$.
\finprf
\begin{coro}
The expansion of $T(k)^{-1}$ as a series in $a$'s takes the form
\be
T_{\infty}^{-1}=\II+\sum_{n=1}^\infty \frac{(-1)^{n}}{n!}
a^\dag_{n\ldots1} \bar T^{(n)}_{\infty1\ldots n}a_{1\ldots n}
\ee
where $\bar T^{(n)}_{\infty1\ldots n}$ is defined by the following inductive expressions:
\bea
\bar T^{(1)}_{\infty0} &=& \II-R_{0\infty}\\
\bar T^{(n+1)}_{\infty01\ldots n} &=& \frac{1}{n+1}
\sum_{i=0}^{n} (\cR^{01\ldots n}_{p_{i}})^{-1}\,
\bar T^{(n)}_{\infty2\ldots i,0,i+1,..,n}\,\cR^{01\ldots n}_{p_{i}}\nonu
&&-\frac{1}{(n+1)!}
\sum_{\sigma\in\cs_{n+1}}(\cR^{01\ldots n}_{p_{n}\circ\sigma})^{-1}\,
\bar T^{(n)}_{\infty\sigma(1)\ldots 
\sigma(n)}R_{\sigma(0)\infty} \,
\cR^{01\ldots n}_{p_{n}\circ\sigma}
\eea
It obeys the corollary \ref{coroformT}, with $R_{\infty\mu(i)}$ replaced by $R_{\mu(i)\infty}$
\end{coro}
\prf
Simple calculation {from} property \ref{adjoint},
theorem \ref{theoTa} and property \ref{Tinv}.
\finprf
\begin{prop}
        The vertex operators $T$ defined in theorem \ref{theoTa} induce an 
        isomorphism between the algebras $\ca_{R}$ and $\ca_{R^{-1}}$. The 
        isomorphism is given by
        \be
        \tau:\ \left\{\begin{array}{lll}
        \ca_{R} & \rightarrow & \ca_{R^{-1}} \\
        a & \mapsto & \hat{a}=T^{-1}a \\
        a^\dag & \mapsto & \hat{a}^\dag=a^\dag T 
        \end{array}\right.
        \ee
\end{prop}
\prf
We first show that $\hat a$ and $\hat a^\dag$ obey the exchange relations 
of $\ca_{R^{-1}}$. We remind that $R_{12}^{-1}=R_{21}$.
\beano
\hat{a}_{1}\hat{a}_{2} &=& T_{1}^{-1}a_{1}T_{2}^{-1}a_{2}
=T_{1}^{-1}T_{2}^{-1}R_{12}a_{1}a_{2} 
= R_{12}T_{2}^{-1}T_{1}^{-1}a_{1}a_{2}\\
&=& R_{12}T_{2}^{-1}T_{1}^{-1}R_{21}a_{2}a_{1} 
= R_{12}T_{2}^{-1}a_{2}T_{1}^{-1}a_{1}= R_{12}\hat{a}_{2}\hat{a}_{1}
\eeano
One does a similar calculation with $\hat{a}_{1}^\dag\hat{a}_{2}^\dag$. 
In the same way, one computes
\beano
\hat{a}_{1}\hat{a}^\dag_{2} &=& T_{1}^{-1}a_{1}a^\dag_{2}T_{2}
=T_{1}^{-1}a^\dag_{2}R_{12}a_{1}T_{2}+T_{1}^{-1}\delta_{12}T_{2}
= T_{1}^{-1}a^\dag_{2}T_{2}a_{1}+\delta_{12}\\
&=& a^\dag_{2}T_{1}^{-1}R_{21}T_{2}a_{1}+\delta_{12}
= a^\dag_{2}T_{2}R_{21}T_{1}^{-1}a_{1}+\delta_{12}
=\hat{a}^\dag_{2}R_{21}\hat{a}_{1}+\delta_{12}
\eeano
This shows that $\ca_{R}$ is embedded into $\ca_{R^{-1}}$. Performing 
the same calculation starting {from} $\ca_{R^{-1}}$ proves that 
$\ca_{R^{-1}}$ is embedded into $\ca_{R}$. There is thus equality of 
the two algebras.
\finprf
\subsubsection*{Reduction to the finite dimensional case}
The above results can be applied to the case without spectral 
parameter. We have to start with a finite dimensional $R$-matrix 
obeying the Yang-Baxter equation 
\be
R_{12}R_{13}R_{23}=R_{23}R_{13}R_{12}
\ee
and a unitarity condition $R_{12}R_{21}=\II$ where, for this section 
only, the spectral parameters are not present. 
The deformed oscillators algebra is then finite dimensional,
 and all the properties stated above are still valid, the proofs 
 following the same lines, 
omitting the integration over the spectral parameters.

Note however that the unitarity condition has still to be fulfilled, 
and this requirement excludes for instance the (triangular) $R$-matrix 
of the finite-dimensional quantum group $\cu_{q}(sl_{2})$.

\sect{Application to integrable systems\label{s:hier}}
\begin{prop}[Hierarchy associated to $\ca_{R}$]\hfill\\
        Let $H^{(n)}$ be defined by
        \be
        H^{(n)}=\int_{-\infty}^\infty dk\ k^n\, a^\dag(k) a(k),\ \forall n=0,1,2,\ldots
        \ee
        $H^{(n)}$ form an Abelian algebra, which defined a hierarchy for the 
        algebra $\ca_{R}$.
        
        The evolution of the $a$ and $a^\dag$ operators under the flow 
        $H^{(n)}$ is given by
        \bea
        e^{itH^{(n)}}\, a(k)\, e^{-itH^{(n)}} &=& e^{-itk^{n}} a(k)\\
        e^{itH^{(n)}}\, a^\dag(k)\, e^{-itH^{(n)}} &=& e^{itk^{n}} a^\dag(k)
        \eea
\end{prop}
\prf
Direct calculation. For instance:
\beano
a_{1}H^{(n)} &=& a_{1} k_{2}^n a_{2}^\dag a_{2}=
k_{2}^n (a_{2}^\dag R_{12} a_{1}a_{2} +\delta_{12}a_{2})
= k_{2}^n a_{2}^\dag R_{12} R_{21}a_{2}a_{1}+ k_{1}^n a_{1}\\
&=& H^{(n)} a_{1} + k_{1}^n a_{1}
\eeano
and thus
${[H^{(n)},a_{1}]}= -k_{1}^n a_{1}$.
\finprf
\begin{prop}
Any well-bred operator $L$ is an integral of motion for the hierarchy:
        \be
        {[L, H^{(n)}]}=0,\ \forall n=0,1,2,\ldots
        \ee
In particular, it is the case of the well-bred vertex operators defined in theorem 
\ref{theoTa}, and the quantum group $\cu_{R}$ generates an infinite 
dimensional symmetry algebra for the hierarchy.
\end{prop}
\prf
$
L_{1}H^{(m)}_{2} = k_{2}^m\ L_{1}a^\dag_{2}a_{2}
=k_{2}^m\ a^\dag_{2}R_{12}L_{1}a_{2} 
= k_{2}^m\ a^\dag_{2}R_{12}R_{21}a_{2}L_{1}
= H^{(m)}_{2}L_{1}
$
\finprf
\begin{rmk}  \rm {From} the example of section \ref{NLS} (see below), and 
which was studied in \cite{qNLS,MRSZ}, we conjecture that to each 
$\ca_{R}$-hierarchy  
corresponds an integrable system already studied in the literature. 
The $a^\dag$ operators in the Fock space representation, in this context, 
will correspond to 
asymptotic states of the system. The correlation functions of the 
system would then be computed using the $a^\dag$ operators.
\end{rmk}

\sect{Fock space and evaluation representations\label{s:fock}}
Associated to the deformed oscillator algebra $\ca_{R}$ comes the 
notion of Fock space:
\begin{defi}
The Fock space $\cf_{R}$ of the $\ca_{R}$ algebra is the module generated by 
the vacuum $\Omega$ such that
\be
a_{i}(k)\Omega=0\ \ \forall\, i=1,\ldots,N\ \forall\, k
\ee
\end{defi}
Now, since one has constructed a quantum group {from} the $\ca_{R}$ 
algebra, it is natural to look at the representations induced by the 
Fock space:
\begin{prop}
        The Fock space $\cf_{R}$ decomposes under the action of the 
        Hamiltonians $H^{(n)}$ into an infinite sum of tensor product of
        evaluation representations of $\cu_{R}$:
        \be
        \cf_{R}=\oplus_{n=0}^\infty \int dk_{1}\cdots dk_{n}\ 
        \theta(k_{1}\leq k_{2}\leq \ldots\leq k_{n})\,\cv_{n}(k_{1},k_{2},\ldots,k_{n})
        \ee
        where $\theta(k_{1}\leq k_{2}\leq \ldots\leq k_{n})$ indicates that the spectral 
        parameters are ordered.
        
        In particular, the representations $\cv_{n}(k_{1},\ldots,k_{n})$ are 
        of dimension $N^n$, and $T$ acts in these spaces by right-multiplication 
        by $R$.
\end{prop}
\prf
Since the Hamiltonians $H^{(n)}$ form a commuting subalgebra of 
$\ca_{R}$, we can consider them as a Cartan subalgebra, and decompose 
$\cf_{R}$ into Cartan-eigenspaces $\cv_{n}(h_{1},h_{2},\ldots)$, 
where $n$ denotes the eigenvalue under $H^{(0)}$ (which turn to be 
still the particle number although we are in the deformed case) and $h_{p}$ 
is the eigenvalue of $H^{(p)}$ ($p>0$). 
Now, since $\cu_{R}$ commutes with these 
Hamiltonians, the eigenspaces are stable under the action of $\cu_{R}$ 
and thus are representations of $\cu_{R}$. 

The vectors in $\cf_{R}$ are linear combinations of monomials 
$a^\dag_{\alpha_{1}}(k_{1})\cdots a^\dag_{\alpha_{m}}(k_{m})\Omega$, 
$\forall m$. 
On the eigenspace $\cv_{n}(h_{0},h_{1},h_{2},\ldots)$, one must consider only 
monomials with $m=n$: this provides only a finite number of terms, 
and the eigenspace is of finite dimension. Moreover, the eigenvalues 
under the $H^{(n)}$ being fixed, one has equations
\beano
h_{1} = \sum_{i=1}^n k_{i} \ ;\
h_{2} = \sum_{i=1}^n k_{i}^2 \ ;\ \ldots\ ;\ 
h_{n} = \sum_{i=1}^n k_{i}^n
\eeano
which completely fixes the values of $k_{1},\ldots,k_{n}$ (up to a 
permutation) and also of 
$h_{p}=\sum_{i=1}^nk_{i}^p,\ p>n$. Thus, we can replace the labeling 
$h_{1},h_{2},\ldots$ by $k_{1},\ldots,k_{n}$, whence the 
notation $\cv_{n}(k_{1},\ldots,k_{n})$ for the representations of 
$\cu_{R}$. Finally, the exchange relations among the $a^\dag$'s allow 
to reorder them in such a way that the spectral parameters are in
increasing order. 

Because it is a vertex operator, the action of $T$ on $\Omega$ is trivial, 
and since it is well-bred, its action on other states is a 
multiplication by $R$.
\finprf

\begin{rmk}[Hopf structure of $\cu_{R}$]\hfill\\
  \rm Although one cannot obtain 
the Hopf structure of $\cu_{R}$ starting {from} $\ca_{R}$, one can 
infer it {from} the present construction in the following way.

The "first" eigenspaces are
\beano
\cv_{0}(0) &=& \CC\,\Omega \\
\cv_{1}(k) &=& \mbox{Span}\Big(a^\dag_{i}(k)\Omega,\ i=1,\ldots,N\Big) \\
\cv_{2}(k_{1},k_{2}) &=& 
\mbox{Span}\Big(a^\dag_{j}(k_{2})a^\dag_{i}(k_{1})\Omega,\
        k_{1}\leq k_{2},\ i,j=1,\ldots,N\Big)
\eeano
Looking at the action of the well-bred vertex operators $T$ on 
these spaces, one gets
\be
T\Omega=\Omega\mb{;}
T_{1}a^\dag_{2}\Omega=a_{2}^\dag R_{12}\Omega\mb{;}
T_{1}a^\dag_{2}a^\dag_{3}\Omega=a_{2}^\dag R_{12}a_{3}^\dag R_{13}\Omega
\ee
Interpreting $\cv_{2}(k_1,k_{2})$ as the tensor product 
$\cv_{1}(k_{1})\otimes \cv_{1}(k_{2})$, 
\be
a_{2}^\dag a_{3}^\dag\Omega\sim a_{2}^\dag \Omega\otimes 
a_{3}^\dag\Omega
\ee
we get\footnote{Be careful that the indices 1, 2, 3 refer to the 
auxiliary spaces while the tensor product refers to $\ca_{R}$.}
\be
T_{1}a^\dag_{2}a^\dag_{3}\Omega=
a_{2}^\dag R_{12}a_{3}^\dag R_{13}\Omega\sim 
a_{2}^\dag R_{12}\Omega\otimes a_{3}^\dag R_{13}\Omega
=(T_{1}\otimes T_{1})\, \big(a_{2}^\dag \Omega\otimes 
a_{3}^\dag\Omega\big)
\ee
Thus,
we are naturally led to the 
coproduct formula
\be
\Delta(T)=T\otimes T
\ee
which is the right one for $\cu_{R}$.
\end{rmk}

\begin{rmk} \rm Note also that, due to the finite number of $a$ 
operators in the states of $\cv_{m}$, the vertex operators truncate at 
level $m$, and turn to be polynomials in $a$, $a^\dag$ in these 
representations.
\end{rmk}

\sect{Examples\label{s:ex}}
We treat here two examples: one associated with an additive spectral 
parameter, and the second one to a multiplicative spectral parameter.
\subsection{The nonlinear Schr\"odinger equation\label{NLS}}
The nonlinear Schr\"odinger equation in 1+1 dimensions (NLS) has been 
widely studied. We look at it in the QISM approach ( for a review, 
see for instance 
\cite{Gut} and ref. therein).

It has already been shown \cite{qNLS, MRSZ} that all the informations on the 
hierarchy associated to NLS can be reconstructed starting {from} the algebra 
$\ca_{R}$, 
where $R$ is the $R$-matrix of the Yangian $Y(N)$ based on $gl(N)$:
\be
R(k)=\frac{1}{k+ig}\left(k\,\II_{N}\otimes \II_{N} 
+ig\, P_{12}\right)\mb{,} P_{12}=\sum_{i,j=1}^N E_{ij}\otimes E_{ji}
\ee
This $R$-matrix obey an additive Yang-Baxter equation
\be
R_{12}(k_1-k_2)R_{13}(k_1-k_3)R_{23}(k_2-k_3)=
R_{23}(k_2-k_3)R_{13}(k_1-k_3)R_{12}(k_1-k_2)
\ee
and one shows, using $P^2=\II$, that $R_{12}(k)R_{21}(-k)=\II$. Thus, 
the properties stated above apply.

In fact, 
it is well-known that
 the canonical field $\Phi$ obeying the (quantum) NLS:
\[
\Big(i\prt_t+\prt^2_x\Big)\Phi(x,t) = 2g\, 
:\Phi(x,t)\bar\Phi(x,t)\Phi(x,t):
\mb{with} \Phi(x,t)=\left(\begin{array}{c} \vph_{1}(x,t) \\ \vdots \\ 
\vph_{n}(x,t) \end{array}\right)
\]
can be reconstructed {from} $\ca_{R}$ \cite{qNLS}. 
The Hamiltonian is then exactly 
$H^{(2)}$, and the Yangian $Y(N)$ is a symmetry of the hierarchy 
\cite{MW,MRSZ}. The operators $a^\dag$ correspond to asymptotic 
states in the Fock space $\cf$. 

The generators $Q_{0}^a$ and $Q_{1}^a$ of $Y(N)$ in its Drinfeld 
presentation were built in term of $\ca_{R}$ in \cite{MRSZ} (see also 
\cite{MW} for the $gl_{2}$ case). The 
present approach is an alternative construction of $Y(N)$ in the FRT 
presentation. It has the advantage to give an explicit construction 
for all the generators of the Yangian, and also to give the action of 
 these generators (\ie of the integrals of motion) on the $a$ and 
 $a^\dag$ operators (\ie the asymptotic states of the system).
 
\subsection{The quantum group $\cu_{q}(\wh{gl_{2}})$}
We take here the evaluated $R$-matrix of the \underline{centerless} affine 
$gl_{2}$ quantum 
algebra. Following the usual notation, the spectral parameter is 
denoted $z$. The $R$-matrix reads:
\begin{equation}
  R(z) = \left( \begin{array}{cccc}
  1 & 0 & 0 & 0 \\
  0 & \displaystyle \frac{q(1-z^2)}{1-q^2z^2} & \displaystyle
  \frac{z(1-q^2)}{1-q^2z^2} & 0 \\
  0 & \displaystyle \frac{z(1-q^2)}{1-q^2z^2} & \displaystyle
  \frac{q(1-z^2)}{1-q^2z^2} & 0 \\
  0 & 0 & 0 & 1 \\
\end{array} \right) \;.
\label{eq:ruq}
\end{equation}
It is defined here up to a normalization factor $\rho$ such that the unitarity condition 
$R_{12}(z_{1}/z_{2})\,R_{21}(z_{2}/z_{1})=1$ is preserved, \ie
\begin{equation}
  \rho(z) \, \rho(\frac{1}{z})=1\;.
\end{equation}
The $R$-matrix obeys a multiplicative Yang-Baxter equation:
\be
R_{12}(z_1/z_2)R_{13}(z_1/z_3)R_{23}(z_2/z_3)=
R_{23}(z_2/z_3)R_{13}(z_1/z_3)R_{12}(z_1/z_2)
\ee
and once again, one can apply the above properties. Note however that 
we are forced to take a vanishing central charge, so that
the algebra $\cu_{q}(\wh{gl_{2}})$ is defined by the relation
\be
  R_{12}(z_1/z_2) \, T_{1}(z_1) \, T_2(z_2)
  =
  T_2(z_2) \, T_1(z_1) \, R_{12}(z_1/z_2) 
\ee
The Hamiltonian $H^{(2)}$ should correspond to the 
Hamiltonian of Sine-Gordon model.

\subsection{The elliptic quantum group $\ca_{q,p}(\wh{gl_{2}})$}
 The elliptic quantum group $\ca_{q,p}(\wh{gl_{2}})_{c}$ has defining 
 relations
 \be
   R_{12}(z_1/z_2;q,p) \, T_{1}(z_1) \, T_2(z_2)
  =
  T_2(z_2) \, T_1(z_1) \, R^{*}_{12}(z_1/z_2;q,p) 
\ee
where $R^{*}_{12}(z;q,p)=R_{12}(z;q,pq^{-2c})$. Note that $R_{12}$ obeys the 
unitarity condition. Thus, in  the 
\underline{centerless} case, one has $R^{*}=R$, and the above 
procedure can be applied. One will start with
 the evaluated $R$-matrix of $\ca_{q,p}(\wh{gl_{2}})_{c=0}$ and 
construct the corresponding ZF algebra. 

In this way, one gets a 
well-bred vertex operator that realizes 
$\ca_{q,p}(\wh{gl_{2}})_{c=0}$, and this latter algebra is a symmetry 
of the hierarchy associated to the ZF algebra. In particular, the 
Hamiltonian $H^{(2)}$ should be related to the XYZ model, and in this 
framework, we naturally gets $\ca_{q,p}(\wh{gl_{2}})_{c=0}$ as a 
symmetry of this model.
 
\sect{Conclusion and perspectives\label{s:concl}}
Starting with any $R$-matrix with spectral parameter, obeying the 
Yang-Baxter equation and a unitarity condition, we have constructed 
 the corresponding quantum group $\cu_{R}$ in term  of a deformed 
oscillators algebra $\ca_{R}$. The 
realization we present is an infinite series, the expansion being 
given in the number of creation operators. Up to a normalization 
constant, the construction is unique.
These "well-bred vertex 
operators" act naturally on $\ca_{R}$. As a consequence, they are 
integrals of motion of the integrable hierarchy naturally associated to 
$\ca_{R}$.

Taking as an example the $R$-matrix of $Y(N)$, the Yangian based on $gl(N)$, 
we recover by this construction the nonlinear Schr\"odinger equation 
and its $Y(N)$ symmetry. It is thus very natural to believe that the 
other integrable systems known in the literature can be treated with 
the present approach.

\null

Of course, the comparison between the vertex operators constructed in 
this paper, and the vertex operators of quantum affine algebras known in the 
literature (e.g. \cite{affOV}) has to be done. Note however that our 
construction can be done for {\em any} infinite quantum group, 
provided its evaluated $R$-matrix obeys the unitarity condition.

\null

As a generalization, it is natural to ask whether such an approach can 
be extended to the case of (elliptic) quantum groups with 
non-vanishing central 
charge: this seems to be very much
the case \cite{AFS}. If such a generalization 
can be done, it would then be possible to look at (off-shell) correlation  
functions for the underlying integrable systems. Moreover, this could 
give a pertinent insight in the research of vertex operators, as they 
are looked for when starting with the canonical fields of the 
integrable system \cite{jap}.

\section*{Acknowledgments}
I would like to thank D. Arnaudon and L. Frappat for fruitful remarks 
on the center of ZF algebras.

I am grateful to the referee for pertinent remarks, specially pointing out 
a mistake in the first version of  lemma \ref{lem.eq} 
and  theorem \ref{theoTa}.



\begin{thebibliography}{99}
\bibitem{ZF} A. B. Zamolodchikov and A. B. Zamolodchikov, Ann. Phys.
{\bf 120} (1979) 253;\\
L. D. Faddeev, Soviet Scientific Reviews Sect. C {\bf 1}
(1980) 107.

\bibitem{MRSZ} M. Mintchev, E. Ragoucy, P. Sorba and Ph. Zaugg, J. Phys.
{\bf A32} (1999) 5885.

\bibitem{MW} S. Murakami and M. Wadati, J. Phys. {\bf A29} (1996) 7903.

\bibitem{qNLS} E. Sklyanin, L. D. Faddeev, Sov. Phys. Dokl.
   {\bf 23} (1978) 902; \\
E. Sklyanin, Sov. Phys. Dokl. {\bf 24} (1979) 107;\\
H.B. Tacker, D. Wilkinson, Phys. Rev. {\bf D19} (1979) 3660;\\
 D.B. Creamer, H.B. Tacker, D. Wilkinson, Phys. Rev. {\bf 
D21} (1980) 1523;\\
 J. Honerkamp, P. Weber, A. Wiesler, Nucl. Phys. {\bf
B152} (1979) 266;\\
 B. Davies, J. Phys. {\bf A14} (1981) 2631.

\bibitem{Gut} E. Gutkin, Phys. Rep. {\bf 167} (1988) 1.

\bibitem{affOV} I. Frenkel and N. Jing, Proc. Nat. Acad. Sci. USA {\bf 
85} (1988) 9373;\\
Bai-Qi Jin, Shan-You Zhou, 
{\it Vertex Operator of $U_q(\widehat{B_l})$ for Level One},
{\tt q-alg/9512005}.

\bibitem{AFS} D. Arnaudon, L. Frappat and E. Ragoucy, work in progress.

\bibitem{jap} Y. Hara, M. Jimbo, H. Konno, S. Odake, J. Shiraishi,
{\it On Lepowsky-Wilson's Z-algebra}, {\tt mathQA/0005203};\\
 Y.Hara, M.Jimbo, H.Konno, S.Odake, J.Shiraishi,
 {\it Free Field Approach to the Dilute $A_L$ Models}, J. Math. 
 Phys. {\bf 40} (1999) 3791, {\tt 
 math.QA/9902150};\\
M. Jimbo, H. Konno, S. Odake, J. Shiraishi,
Comm. Math. Phys. {\bf 199} (1999) 605, {\tt math.QA/9802002}.

\end{thebibliography}
\end{document}